\documentclass[12pt,leqno]{article}
\usepackage{amsmath,amsthm}

\def\init{\setcounter{equation}{0}}
\newtheorem{theorem}{Theorem}[section]
\newcommand{\R}{{\bf R}}

\newtheorem{lemma}{Lemma}[section]

\title{\bf Inverse Boundary Value Problems for Systems of Partial Differential  
Equations.}
\author{Gregory Eskin  and  James Ralston,\\  
 Department of Mathematics, UCLA,\\ Los Angeles,
CA 90095-1555, USA. \ E-mail: eskin@math.ucla.edu}

\begin{document}

\maketitle

\begin{abstract}
We describe the main results and the ideas of the proofs 
in the papers [E] and [ER2] (see References).
In addition, we simplify the construction of asymptotic solutions in [E],
using the results of [ER2], and we simplify the proof of 
estimate (3.9) that was given in [ER2].
\end{abstract}

\section{The Schr\"odinger Equation with an External Yang-Mills Potential.}
\label{section 1}
\init

 Let $\Omega$ be a smooth bounded domain in $\R^n,\ n\geq 3$. Consider the Dirichlet
problem, $u=f$ on $\partial \Omega$, for a system of
differential equations of the form
\begin{equation}                                \label{eq:1.1}
-\Delta u-2i\Sigma_{k=1}^n A_k(x)\frac{\partial u}{\partial x_k} + Bu=0,
\ x\in  \Omega, 
\end{equation}
where $u=(u_1,\dots,u_m)$, and $A_k,\ k=1,\dots,n,$ and $B$ are smooth $m\times m$ matrix
functions. 
We assume that $\Omega$ is such that the Dirichlet problem has a unique solution $u\in    H_1(\Omega)$ for every $f\in H_{\frac{1}{2}}(\partial \Omega)$. 
Let $\Lambda$ be the  Dirichlet-to-Neumann operator,  
\[
\Lambda f=\frac{\partial u}{\partial\hat{n}}+iA\cdot\hat{n} u\mbox{ on }\partial\Omega,
\]
where $\hat{n}$ is the exterior unit normal to 
$\partial\Omega$ and $u$ is the solution of (\ref{eq:1.1}) with 
$u=f$ on $\partial\Omega$.  
The inverse boundary value
problem is to find the coefficients in (\ref{eq:1.1}) given $\Lambda$.
When one rewrites (\ref{eq:1.1}) in the form
\begin{equation}                        \label{eq:1.2}
(-i\frac{\partial}{\partial x}+A(x))^2u +V(x)u=0
\end{equation}
with $A=(A_1,..,A_n)$ and $V=B-\Sigma_{k=1}^n A_k^2+i\Sigma_{k=1}^n \frac{\partial  A_k}{\partial x_k}$, it 
becomes the time-independent Schr\"odinger equation for a particle under the influence of  the   Yang-Mills potential $(A,V)$. We say that two Yang-Mills potentials $(A^{(1)},V^{(1)})$  and   $(A^{(2)},V^{(2)})$ are gauge equivalent if there exists a smooth, invertible matrix  function   $g(x)$ on $\overline{\Omega}$ such that
\begin{equation}                     \label{eq:1.3}
A^{(2)}=g^{-1}A^{(1)}g-ig^{-1}\frac{\partial g}{\partial x}
\mbox{ and }  V^{(2)}=g^{-1}V^{(1)}g.
\end{equation}
The following theorem was proven in [E].
\begin{theorem}                       \label{theo:1.1}
 Let $L^{(j)}u=0,\ j=1,2,$ be two Schr\"odinger equations of the
form (\ref{eq:1.2}) in $\Omega\subset \R^n,\ n\geq 3$, 
with Yang-Mills potentials  $(A^{(j)},V^{(j)}),
\ j=1,2,$ and let $\Lambda^{(j)},\ j=1,2,$
be their Dirichlet-to-Neumann operators. Assume that $\Omega$ is convex.  Then    $\Lambda^{(1)}=\Lambda^{(2)}$ if
and only if $(A^{(1)},V^{(1)})$ and $(A^{(2)},V^{(2)})$
are gauge equivalent.
\end{theorem}
The proof of Theorem \ref{theo:1.1} is based on the method 
of complex exponential solutions with a large
parameter that was introduced in [SU]. 
Let  $\mu$,$\nu$ and $l$ be pairwise orthogonal  vectors 
in ${\R}^n$ with $|\mu|=|\nu|=1$. We set $\theta=\mu+i\nu$, and 
for $\tau >>0$ we set $\zeta =l/2 + (\tau^2-|l|^2/4)^{1/2}\mu$ and 
$\delta=\zeta+i\tau\nu$. We look for solutions of (\ref{eq:1.2}) of the form 
$u=v\exp(ix\cdot \delta)$.
Hence $v$ must satisfy
\begin{equation}                            \label{eq:1.4}
 L_\delta v=_{\mbox{def.}}(-i\frac{\partial}{\partial x} 
+\delta)^2v -2iA(x)\cdot  (\frac{\partial}{\partial x} 
+i\delta)v +B(x)v=
0.
\end{equation}
In order to solve systems of the form (\ref{eq:1.4}) or, more generally,
 inhomogeneous systems
\begin{equation}                             \label{eq:1.5}
L_\delta v =f\mbox{ in }\Omega_0,
\end{equation}
where $\overline{\Omega}\subset \Omega_0$, we will need
solutions of matrix equations of the form
\begin{equation}                           \label{eq:1.6}
i\theta\cdot \frac{\partial C(x,\theta)}{\partial  x}
=\theta\cdot A(x)C(x,\theta),\ x\in\Omega_0,
\end{equation}
where $C$ is an invertible $m\times m$ matrix function. 
One can show that, when $A(x)$ is
extended to a function of compact support in ${\R}^n$, 
there may be no solution of
(\ref{eq:1.6}) which tends to the identity matrix as 
$|x|\rightarrow\infty$. 
However, there are always
solutions which grow polynomially. The following lemma was proven in [E].
\begin{lemma}                              \label{lma:1.1}
 There exists an invertible matrix function $C(x,\theta)$,
  solving (\ref{eq:1.6}) and depending smoothly on 
$(x,\theta)$ on the domain $\overline{\Omega_0}
  \times\{ \theta=\mu+i\nu: \mu\cdot\nu=0,\ |\mu|=|\nu|=1\}$.
This solution is not unique, but it can be chosen to
satisfy $C(x,e^{i\omega}\theta)=C(x,\theta)$.
\end{lemma} 
Using Lemma \ref{lma:1.1}, we can prove (see [E]):
\begin{lemma}                             \label{lma:1.2}
 For any $f\in L^2(\Omega_0)$ there is a $v(\tau)\in H^2(\Omega_0)$  
such that for $\tau$ sufficiently
large $v(\tau)$ solves (\ref{eq:1.5}) in $\Omega$, and
\[
\|v\|_{H^l(\Omega_0)}\leq \frac{C}{ (1+\tau)^{1-l}}\|f\|_{L^2(\Omega_0)},\ l=0,1,
\]
where $C$ is independent of $\tau$ and $f$.
\end{lemma}
To prove Lemma \ref{lma:1.2} we proceed as follows. Let $C(x,\mu+i\nu)$ be as 
in Lemma \ref{lma:1.1}, and let  $c(x,D,\tau)$ 
be the pseudo-differential operator with 
symbol $C(x,\frac{\xi^\prime +\zeta}{  |\xi^\prime +\zeta|} +i\nu)\chi$, 
where $\xi^\prime =\xi-(\xi\cdot\nu)\nu$
and $\chi$ is a suitable cutoff function in $(\zeta,\xi,\tau)$. 
We look for $v$ in the form
\begin{equation}                            \label{eq:1.7}
v=c(x,D,\tau)E(\tau)g,\mbox{ where }E(\tau)=
(-i\frac{\partial}{\partial x}  +\delta)^{-2}\mbox{ and }g\in L^2(\Omega_0).
\end{equation}
Substituting (\ref{eq:1.7}) into (\ref{eq:1.5}) one gets
$g+T(\tau)g=f$, where the norm of $T(\tau)$ goes to zero as 
$\tau\rightarrow\infty$. Therefore
$I+T(\tau)$ is invertible for $\tau$ large. 
Note that the proof of Lemma \ref{lma:1.2} in [E] is a 
generalization of the method in [ER1] which treated the (scalar) case 
of electro-magnetic potentials. Lemma \ref{lma:1.2} is used in [E] 
to prove the following:

\begin{lemma}                         \label{lma:1.3} 
 For every vector of polynomials, $p(z)$, in the complex variable 
$z=\theta\cdot x$ there is a solution $v(\tau)$ of  (\ref{eq:1.4})
satisfying
\begin{equation}                           \label{eq:1.8}
\lim_{\tau
\rightarrow \infty}v(\tau)=C(x,\theta)\Pi^+(C^{-1}p),
\end{equation}
where $C$ is the matrix function from Lemma \ref{lma:1.1}
and $\Pi^+$ is a Toeplitz projection. 
\end{lemma}

Using (\ref{eq:1.7}) to construct the leading order term in 
$v(\tau)$ made both the proof of Lemma \ref{lma:1.3} and its
applications somewhat
complicated in [E]. In [ER2] we found a way to construct 
such solutions more simply with explicit higher order 
asymptotics as $\tau\to \infty$, and this lead to simpler proofs.
 We will give this construction in \S 2. 

Now we can complete the proof of Theorem \ref{theo:1.1}. 
Following the strategy of [SU], we
can use the assumption $\Lambda^{(1)}=\Lambda^{(2)}$, 
Green's formula and 
Lemma \ref{lma:1.3} to 
derive integral identities involving $C^{(j)}$ and $(A^{(j)},V^{(j)}),\ j=1,2$. 
Then
arguments involving $\overline{\partial}$-equations in the   parameters in 
these identities (see \S 5 and \S 6 in [E]) lead to the proof 
of Theorem \ref{theo:1.1}.

\section{Construction of Complex Exponential Solutions.}
\label{section 2}
\init

To construct solutions of (\ref{eq:1.1}) of the form $u=v\exp(ix\cdot \delta)$ 
we will  proceed as follows. 
Substituting $u=v\exp(ix\cdot \delta)$ into (\ref{eq:1.1}), 
one sees that $v$ must satisfy
\begin{equation}                                   \label{eq:2.1}
L_\delta v=_{\mbox{def.}}-\Delta v-2i\delta\cdot\frac{\partial v}{\partial x}
-2i A\cdot(i\delta v +\frac{\partial v}{\partial x})+Bv=0.
\end{equation}
We will construct $v$ in the form $v=\Sigma_{k=0}^nv_k +\tilde{v}_n$, where
$v_k=O(\tau^{-k})$ and is explicit modulo solutions of (\ref{eq:1.6}), and 
$\tilde v_n$ is $O(\tau^{-n-1})$. We have 
\[
\delta =\tau\theta +\frac{l}{2} +O(\tau^{-1})
=_{\mbox{def.}}\tau\theta +\delta^\prime  \mbox{ and }
\]
\[
L_\delta
=_{\mbox{def.}}-2i\tau\theta\cdot\frac{\partial}{\partial x} 
+2\tau \theta\cdot A  +
M_{\delta^\prime}.
\]
To solve (\ref{eq:2.1}) modulo terms of order $\tau^{-n}$ we require
\[
2\tau(i\theta\cdot\partial_x -\theta\cdot A)v_k=M_{\delta^\prime}v_{k-1}
\]
for $k=0,\dots,n$ with $v_{-1}=0$. We set $v_0=C_0(x,\theta)p(x\cdot \theta)$, 
where $C_0$  is a solution of (\ref{eq:1.6}) with the
properties described in Lemma \ref{lma:1.1} and $p(z)$ is a vector of 
polynomials in the complex  variable $z$. Since we only require 
that $v$ satisfy (\ref{eq:2.1}) on $\Omega$, we can introduce 
a  cut-off function $\psi\in C_0^\infty(\Omega_0)$, $\psi =1$ on 
a neighborhood of $\Omega$, and set
\begin{equation}                        \label{eq:2.2}
v_k=-i(2\tau)^{-1}C_k(x,\theta)(\theta\cdot\partial_x)^{-1}
(C_k^{-1}(\cdot,\theta)\psi  M_{\delta\prime}v_{k-1}),\ k=1,\dots,n.
\end{equation}
In (\ref{eq:2.2}) the operator 
$(\theta\cdot\partial_x)^{-1}$ multiplies the Fourier transform
 by  $(i\theta\cdot \xi)^{-1}$, and $C_k$ is again a solution 
of (\ref{eq:1.6}) as in Lemma \ref{lma:1.1} 
(in  applications so far we have taken $C_k=C_0$,
 but this is not necessary). 
Since  $M_{\delta^\prime}$ does not increase the order of a term in $\tau$
and $i(2\tau)^{-1}(\theta\cdot\partial_x)^{-1}$  adds a factor of $\tau^{-1}$, 
we have
$v_k=O(\tau^{-k})$ and 
\[
L_\delta(v_0+v_1+\cdots +v_n)=M_{\delta^\prime }v_n=O(\tau^{-n})
\]
on the neighborhood of $\Omega$ where $\psi=1$. 
Hence, taking $\psi_0\in C^\infty(\Omega_0)$ supported in the set 
where $\psi=1$ such that $\psi_0=1$ 
on $\Omega$, $v=v_0+\cdots +v_n+\tilde v_n$ will be a solution of (\ref{eq:2.2}) 
in $\Omega$ if
\begin{equation}                        \label{eq:2.3}
L_\delta\tilde v_n=-\psi_0 M_{\delta^\prime} v_n
\end{equation}
in $\Omega_0$. 

To solve (\ref{eq:2.3}) apply Lemma \ref{lma:1.2}.
 Lemma \ref{lma:1.2} holds in a more general form (see [ER2]):
if $f\in  H^k(\Omega_0),\ k\geq 0$, then $v\in H^{k+2}(\Omega_0)$, 
and one has the estimate
\[
\|v\|_{H^{k+l}(\Omega_0)}\leq \frac{C}{(1+\tau)^{1-l}}\|f\|_{H^k(\Omega_0)},\ l=0,1.
\]
Since $v_n$ is bounded by $C\tau^{-n}$, the solutions given by 
Lemma \ref{lma:1.2}
 will be bounded by $C\tau^{-n-1}$ in $H^k(\Omega_0)$ for $k\geq 0$.
Thus the $v(\tau)$ that we have constructed is a 
solution of (\ref{eq:2.1}) whose asymptotics in $H^k(\Omega)$ up to order 
$\tau^{-n}$ are given by the  asymptotics of
$v_0+\cdots +v_n$. For the leading term in the asymptotics we have 
\begin{equation}                            \label{eq:2.4}
v(\tau)=C_0(x,\theta)p(\theta\cdot x)+O(\tau^{-1}),
\end{equation}
where $O(\tau^{-1})$ means bounded by $C\tau^{-1}$ in $H^k(\Omega), k\geq 0$. 
Since the
limit of $v(\tau)$ no longer involves a Toeplitz projection, 
one no longer needs one of
the arguments (Lemma 5.1) in [E].

\section{The Equations of Isotropic Elasticity.}
\label{section 3}
\init

The construction presented in \S 2 can be used to construct solutions of
 the system of isotropic elasticity as well. Using subscripts for derivatives, 
this system is given  by 
\begin{equation}                         \label{eq:3.1}
 \Sigma_{j=1}^3(\lambda w^j_{x_j})_{x_k} + \Sigma_{j=1}^3(\mu(w^k_{x_j}
+w^j_{x_k}))_{x_j}=0,\ k=1,2,3,
\end{equation}
where $w=(w^1,w^2,w^3)$ is the deformation of an elastic body with 
\lq\lq Lam\'e  parameters"
$\lambda(x)$ and $\mu(x)$. Let $w$ be the solution of the Dirichlet problem 
for (\ref{eq:3.1}) with $w=h$ on $\partial \Omega$. Then the inverse 
boundary value problem for this system is to
recover $\lambda$ and $\mu$ from the Dirichlet-to-Neumann map
\[
\Lambda(h)^k=\Sigma_{j=1}^3(\lambda   w^j_{x_j})\hat{n}^k+\Sigma_{j=1}^3\mu(w^k_{x_j}+w^j_{x_k})\hat{n}^j,
\  k=1,2,3.
\]
There is as yet no proof that $\Lambda$ determines $\lambda$ and $\mu$. 
Partial results are  given in [NU2] and [ER2]. The system (\ref{eq:3.1}) 
is not in the form (\ref{eq:1.1}). 
However, Ang, Ikehata, Trong  and Yamamoto in [AITY] show 
that $w=\mu^{-1/2}u+\mu^{-1}\nabla f-f\nabla \mu^{-1}$ will
satisfy (\ref{eq:1.1}) when the 4-vector $(u,f)$ satisfies the system
\begin{equation}                               \label{eq:3.2}
\Delta\left(\begin{array}{c} u\\ f\end{array}\right) 
+V_1(x)\left(\begin{array}{c} \nabla f\\ \nabla\cdot  u\end{array}\right)
+V_0(x)\left(\begin{array}{c} u\\ f \end{array}\right)
=\left(\begin{array}{c} 0\\  0\end{array}\right).
\end{equation}
Here
\[
V_1(x) =\left(\begin{array}{ll} -2\mu^{1/2}\nabla^2\mu^{-1}
& -\mu^{-1}\nabla \mu\\
0
& \frac{\lambda +\mu}{\lambda +2\mu}\mu^{1/2}\end{array} \right)
\]
and $\nabla^2f$ denotes the Hessian matrix 
$\partial^2f/\partial x_j\partial x_k$. 
The matrix $V_0$ is a complicated expression in $\lambda$, $\mu$ 
and their derivatives, but  it
vanishes when $\mu$ is constant. Since (\ref{eq:3.2}) does have 
the form (\ref{eq:1.1}), the method 
above can be used to construct solutions $(u,f)=\exp(i\delta\cdot x)(r,s)$ 
of (\ref{eq:3.1}) with  prescribed 
asymptotics as $\tau\rightarrow\infty$.

Suppose that one has two elastic bodies occupying the region $\Omega$ 
with Lam\'e parameters $(\lambda^{(j)},\mu^{(j)}).\ j=1,2$. 
Let $\Lambda^{(j)},\ j=1,2,$ be the
corresponding Dirichlet-to-Neumann maps. Using Green's formula, 
one can verify that
$\Lambda^{(1)}=\Lambda^{(2)}$ is equivalent to 
\begin{eqnarray}                             \label{eq:3.3}
  H(w^{(2)},w^{(1)})
=_{\mbox{def.}}\int_\Omega[(\lambda^{(2)}-\lambda^{(1)})
(\nabla \cdot\overline{w^{(2)}})(\nabla \cdot w^{(1)})
\\
+ \frac{1}{2}(\mu^{(2)}-\mu^{(1)})
\Sigma_{1\leq j,k\leq 3}
(\overline{w^{(2),j}_k+w^{(2),k}_j})(w^{(1),j}_k+w^{(1),k}_j)]dx=0
\nonumber
\end{eqnarray}
for all solutions to $L^{(1)}w^{(1)}=0,\ L^{(2)}w^{(2)}=0$, where now 
$L^{(j)}w^{(j)}=0$
is the system (\ref{eq:3.1}) with 
$\lambda=\lambda_j$ and $\mu=\mu_j$. We can use (\ref{eq:3.3})
in the following way.
The method of \S 2 yields solutions of (\ref{eq:3.2}) 
with
 prescribed asymptotics. 
For the
problem corresponding to $L^{(1)}$ we use these with $\delta^{(1)}=\delta$ as
defined earlier, but for the problem corresponding to $L^{(2)}$ we set
$\delta^{(2)}=\overline{\delta} -l$. Taking $w^{(j)}=
\mu_j^{-1/2}u^{(j)}+\mu_j^{-1}\nabla f^{(j)}-f^{(j)}\nabla \mu_j^{-1}$,
 we substitute 
$w^{(1)}$ and $w^{(2)}$ into (\ref{eq:3.3}). 
When we collect the terms of each order in $\tau$, this 
gives
\[
0=H(w^{(2)},w^{(1)})=\tau^2H_2+\tau^1H_1+H_0+ \cdots + \tau^{-N}H_{-N}+O(\tau^{-N-1}), 
\]
where each $H_j$ is independent of $\tau$. How far one can continue this
expansion depends on the choice of $n$ in \S 2. Note that each $H_j$ must vanish when
$\Lambda^{(1)}=\Lambda^{(2)}$. In particular, one has
\begin{equation}                             \label{eq:3.4}
H_2=\int_\Omega e^{il\cdot x}(\theta\cdot\overline{r}_0^{(2)},
\overline{s}_0^{(2)})V(x,\theta)(\theta\cdot r_0^{(1)},s_0^{(1)})^tdx=0,
\end{equation}
where
\[
V(x,\theta)=\left(\begin{array}{ll}
(\lambda_1+\mu_1-\lambda_2-\mu_2)
\frac{(\mu_1\mu_2)^{1/2}}{(\lambda_1+2\mu_1)(\lambda_2+2\mu_2)}
&   2(\mu_2^{-1}-\mu_1^{-1})\mu_2^{-1/2}\theta\cdot\partial_xb_2
\\ 
2(\mu_2^{-1}-\mu_1^{-1})\mu_1^{-1/2}\theta\cdot\partial_xb_1
&  2(\mu_2^{-1}-\mu_1^{-1})(b_1a_1+b_2a_2)\end{array}\right).
\]
The functions $a_j$ and $b_j,\ j=1,2$ are given by
\[
a_j=(\theta\cdot\partial_x)^2\mu_j^{-1}\ \mbox{ and }\ b_j=\frac{\mu_j}
{ 2}\frac{\lambda_j+\mu_j}{\lambda_j+2\mu_j}.
\]
The functions $(r_0^{(j)},s_0^{(j)})$ are vector solutions of 
the following version of (\ref{eq:1.6}) 
\begin{equation}                     \label{eq:3.5}
-2\theta\cdot \partial_x\left(
\begin{array}{c}r_0\\s_0\end{array}\right)=V_1\left(\begin{array}{c} s_0\theta^t\\
\theta\cdot r_0\end{array}\right).
\end{equation}
This is (\ref{eq:1.6}) with   
\[ 
A=-\frac{i}{2}V_1\left(\begin{array}{ll}0_{3\times 3}&\theta^t\\
\theta&0 \end{array}\right).
\] 
One can analyze $H_2$
by the techniques that were used for the Yang-Mills case in [E]. However, here that 
does not lead to the conclusion that $\lambda_1=\lambda_2$ and $\mu_1=\mu_2$. Instead
one arrives at the identity (Theorem 1.3 of [ER2])
\[
b_1^{1/2}(\theta\cdot\partial_x)^2(b_1^{-1/2})-b_1(\theta\cdot\partial_x)^2
 (\mu_1^{-1})=b_2^{1/2}(\theta\cdot\partial_x)^2(b_2^{-1/2})
-b_2(\theta\cdot\partial_x)^2
 (\mu_2^{-1}).
\]
Since this holds identically in $\theta$, it is equivalent to a system of five 
partial differential equations satisfied by $(\lambda_1,\mu_1,\lambda_2,\mu_2)$.
However, even in the case that $\lambda_1$ and $\mu_1$ are constant, these equations
have solutions with $\lambda_2$ and $\mu_2$ nonconstant. 

A more readily useful result is the following (Theorem 2 of [ER2]):
\begin{theorem}                                \label{theo:3.1}
 Let 
\[
\|f\|_\alpha^2=\int_\Omega |f|^2e^{2\alpha
|x|^2}dx.
\]
Then, if $\Lambda^{(1)}=\Lambda^{(2)}$, for $\alpha>\alpha_0$, one has
\[
\|\lambda_1+\mu_1-\lambda_2-\mu_2\|_\alpha\leq 
\frac{C}{\alpha}\|\mu_1-\mu_2\|_\alpha
\]
with $C$ independent of $\alpha$. 
\end{theorem}
This result comes from $\overline{\partial}$-equations 
that arise when one follows the method of \S 6 in [E]. 
Theorem \ref{theo:3.1} can be used to deduce uniqueness for constrained forms of
the inverse problem. For instance, if one is given
either that $\lambda_1=\lambda_2$ or that $\mu_1=\mu_2$, 
and that the Dirichlet-to-Neumann maps are equal, then both Lam\'e parameters 
must be equal.

The second way that one can use (\ref{eq:3.4}) is to choose 
$\theta$ as a function of $l/|l|$ keeping
$\theta(l/|l|)\cdot l=0$. This makes
(\ref{eq:3.4}) 
(and each of the other equations in the sequence 
$\{H_j=0\}_{j=2}^{-\infty}$) equivalent to a pseudo-differential  
equation of the form
\begin{equation}                           \label{eq:3.6}
P(x,D)(\lambda_2-\lambda_1)+Q(x,D)(\mu_2-\mu_1)=0
\end{equation}
To prove uniqueness for the inverse boundary value problem 
we would like to use
(\ref{eq:3.6}) to bound $\lambda_2-\lambda_1$ in terms of $\mu_2-\mu_1$
 or to bound $\mu_2-\mu_1$ in terms of $\lambda_2-\lambda_1$.
Boundary determination for this inverse problem (see [NU1]) 
implies that $\lambda_2-\lambda_1$ and $\mu_2-\mu_1$ vanish to infinite
 order on $\partial \Omega$.
Nonetheless, (\ref{eq:3.6}) does not imply any estimates of this kind 
without more information
on the operators $P$ and $Q$.
Note that the construction in \S 2 makes $(r,s)$ a function of $\theta$, so that
it contributes to the symbols of $P$ and $Q$,and, since the solutions of 
(\ref{eq:3.5}) are not 
explicit in general, one usually does not know what $P$ and $Q$ are.
However, (\ref{eq:3.6}) becomes useful when one assumes that 
$\nabla \mu_j$ is small
in $C^\infty$. In that case all entries in $V_1$ except $(V_1)_{44}$ 
become small and
one can solve (\ref{eq:3.5}) explicitly modulo small terms. 
In fact (\ref{eq:3.5}) has a unique matrix solution tending to
the identity as $|x|\rightarrow \infty$ and given by
\begin{equation}                               \label{eq:3.7}
C_{00}(x,\theta)=
\left(\begin{array}{ll} I_{3\times 3}
&  0  \\
\varphi\theta
&   1     \end{array} \right),
\end{equation}
modulo small terms,
where $\varphi=(\theta\cdot \partial_x)^{-1}(-\psi b\mu^{-1})$.
Choosing 
\[
p=(\hbox{Re}\{\theta(l/|l|)\},0)
\]
 in (\ref{eq:2.4}),
one gets
\begin{equation}                              \label{eq:3.8}
(r^{(j)},s^{(j)})=(\hbox{Re}\{\theta(l/|l|)\},0)+
(\tilde{r}_0^{(j)},\tilde{s}_0^{(j)})+ O(\tau^{-1}),
\end{equation}
where $\tilde{r}_0^{(j)},\tilde{s}_0^{(j)}$
are symbols of order zero and $\tilde{r}_0^{(j)}$ is small.
Since $\theta\cdot\hbox{Re}\{\theta\}=1$, (\ref{eq:3.8}) makes $P(x,D)$ in (\ref{eq:3.6}) 
simply multiplication by 
$(\mu_1\mu_2)^{1/2}(\lambda_1+2\mu_1)^{-1}(\lambda_2+2\mu_2)^{-1}$
plus an operator of order zero whose norm can be made 
arbitrarily small by taking $\nabla \mu_j,\ j=1,2,$ sufficiently
small in $C^k(\Omega)$. Thus (\ref{eq:3.6}) implies
\begin{equation}                          \label{eq:3.9}
\|\lambda_2-\lambda_1\|_{H^k(\Omega)}\leq C_k\|\mu_2-\mu_1\|_{H^k(\Omega)}
\end{equation}
for all $k$, when $\nabla \mu_j,\ j=1,2,$ is sufficiently small.

Uniqueness for the inverse boundary value problem when $\nabla \mu_j$ is small will
follow from (\ref{eq:3.9}) if we can find another estimate bounding a norm of
 $\mu_2-\mu_1$ by a small constant times a norm of $\lambda_2-\lambda_1$. 
One can get this estimate by using $H_0=0$ in the following way. 
The contributions to $H_0$ come from the
expansion of $(r,s)$ in $\tau$ up to order $\tau^{-2}$, and are therefore quite
complicated in general. Nonetheless, 
using (\ref{eq:3.7}) and taking $p=(0,0,0,1)$  in (\ref{eq:2.4}), 
we find that (\ref{eq:3.5}) 
has the solution $(r_0,s_0)=(0,0,0,1) +(\tilde{r}_0,\tilde{s}_0)$, 
where $(\tilde{r}_0,\tilde{s}_0)$ and its derivatives can be made 
arbitrarily small by taking 
$\nabla \mu_j,\ j=1,2,$ sufficiently small in $C^k$-norm. 
Moreover, checking further, one sees that
the higher order terms in the expansion of $(r,s)$ that 
one constructs using (\ref{eq:2.2}) are
also small under these hypotheses. When one uses this choice of $(r,s)$ 
in the
construction of $w^{(1)}$ and $w^{(2)}$, one gets
\begin{eqnarray}  \nonumber                       
H_0=
\int_\Omega e^{il\cdot x}[(\lambda_1-\lambda_2)A(x,\theta(l/|l|),l)\\+
(\mu_1-\mu_2)[(2\mu_1\mu_2)^{-1}|l|^4
+ B(x,\theta(l/|l|),l)]]dx,\nonumber
\end{eqnarray}
where $A$ and $B$ are small symbols in $l$ of orders two and four 
respectively. Thus we have (\ref{eq:3.6}) with 
$\|P(x,D)(\lambda_2-\lambda_1)\|_{L^2(\Omega)}$ bounded by a small 
constant times $\|\lambda_2-\lambda_1\|_{H^2(\Omega)}$ and 
\[
Q(x,D)(\mu_2-\mu_1)=(2\mu_1\mu_2)^{-1}(\Delta)^2(\mu_2-\mu_1)
+R(x,D)(\mu_2-\mu_1),
\]
where $\|R(x,D)(\mu_2-\mu_1)\|_{L^2(\Omega)}$ is bounded 
by a small constant times
$\|\mu_2-\mu_1\|_{H^4(\Omega)}$. Thus we have 
\begin{equation}                           \label{eq:3.10}
\|\mu_2-\mu_1\|_{H^4(\Omega)}\leq\varepsilon (\|\mu_2-\mu_1\|_{H^4(\Omega}
+\|\lambda_1-\lambda_2\|_{H^2(\Omega)}),
\end{equation}
where $\varepsilon$ can be taken arbitrarily small when 
$\nabla \mu_j,\ j=1,2$ is sufficiently small in $C^k(\Omega)$. 
Combining (\ref{eq:3.9}) and (\ref{eq:3.10}), we have the following 
result.

\begin{theorem}                           \label{theo:3.2} 
Given that $\lambda_j,\mu_j$ and $\mu_j^{-1},\ j=1,2$, 
belong to a bounded set, $B$,
in $C^k(\Omega)$ for $k$ sufficiently large, there is an $\varepsilon(B)>0$ such
that $||\nabla \mu_j||_{C^{k-1}(\Omega)}<\varepsilon(B),\ j=1,2,$ implies
$(\lambda_1,\mu_1)=(\lambda_2,\mu_2)$, if $\Lambda^{(1)}=\Lambda^{(2)}$. 
\end{theorem}
This is  the main result of [NU2] and it is also  Theorem 1 of [ER2].
The first version of [ER2] (July 2001) contains Theorem \ref{theo:3.2} with 
additional hypothesis
$\|\nabla \lambda_j\|_{C^{k-1}(\Omega)}<\varepsilon(B),\ j=1,2.$
We did not notice that our proof did not use this hypothesis
until we received a preprint of [NU2] (November 2001).
We are grateful to G.Nakamura and G.Uhlmann for sending this to us.

\end{document}